\documentclass[leqno,12pt]{article}
\usepackage{latexsym}
\usepackage{amssymb}
\usepackage{amsmath}
\usepackage{color}
\usepackage{bbold}
 \usepackage{graphicx}
 \usepackage{subfigure}

 0

\newcommand{\bea}{\begin{eqnarray*}}
\newcommand{\eea}{\end{eqnarray*}}
\newcommand{\be}{\begin{eqnarray}}
\newcommand{\ee}{\end{eqnarray}}
\newcommand{\beq}{\begin{equation}}
\newcommand{\eeq}{\end{equation}}
\newcommand{\dn}{\stackrel{\cal D}{\longrightarrow}}
\newcommand{\ds}{\stackrel{\cal D}{=}}
\newcommand{\Z}{\mathbb{Z}}
\newcommand{\N}{\mathbb{N}}
\newcommand{\R}{\mathbb{R}}
\newcommand{\Cov}{\operatorname{Cov}}
\newcommand{\E}{\mathbb{E}}
\newcommand{\Var}{\operatorname{Var}}

\newtheorem{satz}{Theorem}[section]

\newtheorem{prop}[satz]{Proposition}

\newtheorem{theorem}[satz]{Theorem}

\newtheorem{corollary}[satz]{Corollary}

\DeclareMathOperator*{\argmin}{argmin}

\usepackage[authoryear]{natbib}
\begin{document}

\title{Multiscale change point detection for  dependent data}
\author{
{\small Holger Dette, Theresa Sch\"uler} \\
{\small Fakult\"at f\"ur Mathematik} \\
{\small Ruhr-Universit\"at Bochum} \\
{\small 44799 Bochum, Germany} \\
\and
{\small Mathias Vetter} \\
{\small Mathematisches Seminar} \\
{\small Christian-Albrechts-Universit\"at zu Kiel} \\
{\small 24098 Kiel, Germany} \\
}

\maketitle
\begin{abstract}
In this paper we study the  theoretical properties of  the  simultaneous multiscale change point estimator (SMUCE) proposed  by \cite{frick2014}
in regression models with dependent error processes. Empirical studies show  that in this case the change point estimate is inconsistent,
but it is not known if alternatives  suggested in the literature for correlated data are consistent.  We  propose 
a modification of SMUCE  scaling the basic  statistic  by the long run variance of the error process, which is estimated by 
a difference-type variance estimator calculated from local means from different blocks. For this modification we
 prove model  consistency for physical dependent error  processes and illustrate the finite sample performance by means of a simulation study.
\end{abstract}


\section{Introduction} 
\def\theequation{1.\arabic{equation}}
\setcounter{equation}{0}

The problem  of detecting multiple abrupt changes in the structural properties of a time series 
and to split the data  into several ``stationary''  segments has been of
interest to statisticians for many decades. An efficient a posteriori change-point detection rule
enables the researcher to analyze  data under the assumption of  piecewise-stationarity and has numerous applications 
 including bioinformatics,  neuroscience, genetics,  the  analysis of  speech signals, financial,  and climate data.
 Because of its  importance the literature on the subject is very vast and we refer exemplarily to the work of \cite{yao1988}, \cite{bai1998,bai2003}, \cite{braun2000}, \cite{lavielle2000}, \cite{kolaczyk2005}, 
\cite{davis2006}, \cite{harchaoui2010},
\cite{ciuperca2011,ciuperca2014}, \cite{killick2012}, \cite{fryzlewicz2014}, \cite{matteson2014}, \cite{cho2015}, \cite{preuss2015}, \cite{yau2016}, \cite{haynes2017}, \cite{korkas2017} and \cite{chakar2017}.  This list of references is by no means complete and   further references can be found in the cited literature. 

The focus of  the present  paper is on  the  simultaneous multiscale change point estimator (SMUCE), which was introduced recently in a seminal paper
  of \cite{frick2014} to identify  multiple changes 
 in the mean structure of  the  sequence 
\begin{align} \label{model_1}
Y_i=\vartheta^\ast\left(\tfrac{i}{n}\right)+\varepsilon_i, \quad i=1,\ldots,n,
\end{align}
where $\vartheta^\ast :  [0,1] \to \R$  is  a piecewise constant function and 
$ \varepsilon_{1}, \ldots ,  \varepsilon_{n}$  are independent identically distributed centered Gaussian random variables. Note that these authors considered distributions from a one-parametric exponential
 family with a piecewise constant parameter $\vartheta^\ast$, but for the sake of brevity we restrict ourselves to the location scale model, which corresponds to the Gaussian case.  
 The SMUCE procedure controls  the probability of overestimating the true number of change points, and  it is also possible to give bounds
 for the probability of underestimation. Moreover, one can construct asymptotic honest confidence sets for the unknown step function $\vartheta^\ast$
 and its change points. The method has turned out to be very successful
 and has therefore been extended  in various directions. For example,  \cite{pein2017} consider model \eqref{model_1}
 with a heteroscedastic Gaussian noise process. \cite{li2016} argue that in 
situations with low signal to noise ratio  or with many change-points
compared to the number of observations  SMUCE necessarily leads to a conservative estimate and propose to control 
 the false discovery instead of the family wise error rate. More recently \cite{li2018} extend the procedure 
 to certain function classes beyond step functions in a nonparametric regression setting.

\begin{figure} [ht]
  \centering
\subfigure{
\includegraphics[width=0.45\textwidth, height=160px]{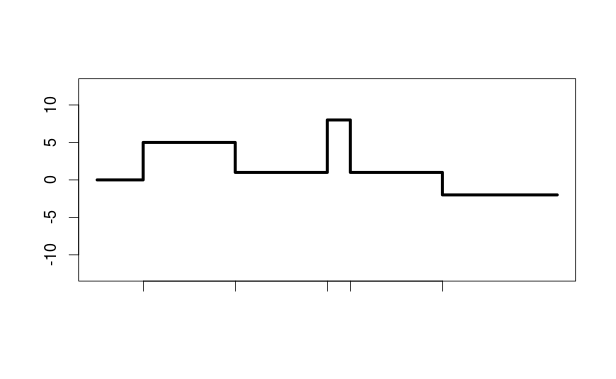}} \quad
\subfigure{
\includegraphics[width=0.45\textwidth, height=160px]{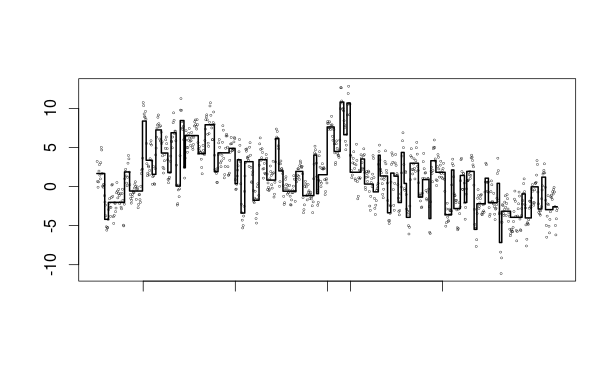}} \\
\vspace{-1.711cm}
\subfigure{
\includegraphics[width=0.45\textwidth, height=160px]{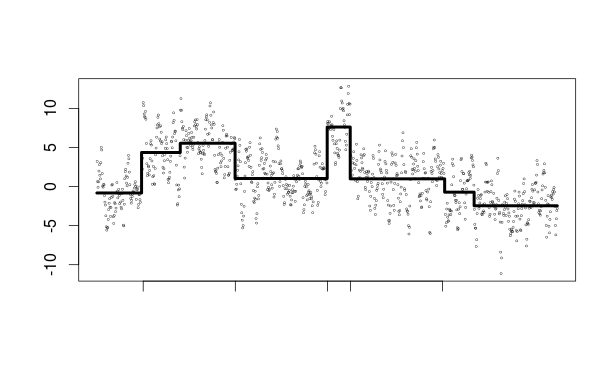}} \quad
\subfigure{
\includegraphics[width=0.45\textwidth, height=160px]{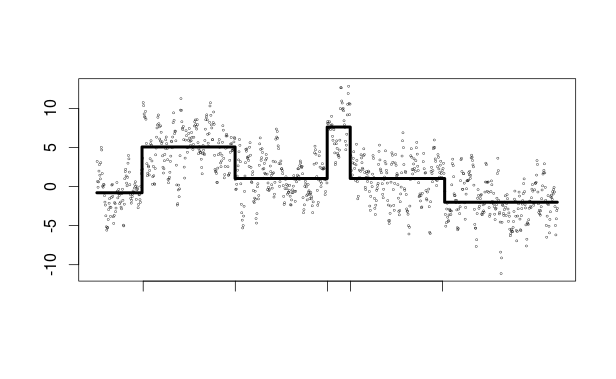}} \quad
\vspace{-1.11cm}
  \caption{\label{Fig1} \it 
 Different estimates of a piecewise constant signal in model \eqref{model_1} with an ARMA$(2,6)$ error process. Upper left panel: true function. Upper right panel: SMUCE. Lower left panel:  estimate proposed in
  \cite{tecuapetla2017}. Lower right panel: estimate proposed in this paper.}
\end{figure}

The present paper is devoted to the analysis of SMUCE in the location scale model \eqref{model_1}  with a piecewise constant regression
function  under more general assumptions on the error process. We are particularly interested in the situation
where the errors  are neither Gaussian nor independent.
If the sample size is reasonably large and the errors are independent, SMUCE is relatively robust because it is based on local means which are asymptotically Gaussian due to the CLT. However, the independence of the errors is more crucial and ignoring this  assumption may lead to serious errors in the estimation procedure.
This is illustrated in Figure \ref{Fig1}, where we display a typical  estimate of the  signal  (upper left panel) by the modification of SMUCE proposed in
 \cite{tecuapetla2017} for $m$-dependent errors (lower left panel). The data generating process is an ARMA$(2,6)$ process.
 We observe that the modification still produces  a function with too many jumps.
 The lower right panel shows the estimate proposed in this paper, which seems to work better. The upper right panel shows the performance of SMUCE, which clearly overestimates the true number of change points.
  A more detailed comparison will be 
presented  in Section \ref{sec4}.
 
The reason for the differences consists in the fact that in the case of dependent data all described procedures require a reliable estimate of the long run 
variance of the error distribution.  \cite{tecuapetla2017}  demonstrate by means of a simulation study  that  the problem can easily be addressed for $m$-dependent errors
using difference based estimators [see \cite{hall1990} or \cite{dette1998}]. Their approach provides a solution for a specific 
error structure and we see the improvement in Figure \ref{Fig1}.
However, from a practical point of view the method requires  a good choice of $m$, and the example  indicates that this procedure might not work well 
for  other dependence structures. More importantly, from a theoretical point of view rigorous statements  regarding  the performance of SMUCE in models with more general (stationary) error processes are missing. It turns out that results of this type are substantially more
difficult to obtain and are--to our best knowledge--not available in the literature so far. 

In this paper we address this problem and prove consistency of SMUCE with an appropriately modified variance estimator
 under the  assumption that the error process  $\{ \varepsilon_i\}_{i\in\Z}$ is a physical system in the sense of 
\cite{wu2005}. This  includes such  important examples as ARMA or  GARCH processes. We also avoid any distributional assumptions regarding the errors
$\varepsilon_i$ except the existence of moments.
In Section \ref{sec2} we introduce the model and the modification of the SMUCE procedure to address general time dependent error processes.
 Roughly speaking, we have to define consistent estimates of the long run variance
\begin{align} \label{lrv}
\sigma_\star^2:=\sum_{k\in\Z} \Cov(\varepsilon_0,\varepsilon_k)~,
\end{align}
which address the fact that the regression function may be only  piecewise  constant and not constant. This is achieved by a two step estimator
which is defined as a  difference based  estimator of local averages. The asymptotic properties of the modified procedure are established in Section \ref{sec3}. We prove that
 the number of change points is identified with probability converging
to $1$ and that all change points are estimated consistently. The finite sample properties are investigated in Section \ref{sec4} by means of a simulation study. 
 Finally, all proofs and technical details are deferred to an appendix.

\section{Multiscale change point detection for dependent data}
 \label{sec2}
\def\theequation{2.\arabic{equation}}
\setcounter{equation}{0}

We begin with a brief review of the simultaneous multiscale change point estimator (SMUCE) as introduced by   \cite{frick2014}, where
we directly  address the problem of  dependent data.  Throughout this paper let  
\begin{align} \label{21} 
\vartheta^\ast(t):=\sum_{k=0}^{K^\ast} \theta_k^\ast \mathbb{1}_{[\tau_k^\ast,\tau_{k+1}^\ast)}(t)
\end{align}
denote the ``true'' unknown signal in model \eqref{model_1}, 
where $K^\ast$ is the (unknown) number of change points, $0=\tau_0^\ast<\tau_1^\ast<\ldots<\tau_{K^\ast}^\ast<\tau_{K^\ast+1}^{\ast}=1$ are the change point locations, and $\theta_0^\ast,\ldots,\theta_{K^{\ast}}^\ast$ are the function values of $\vartheta^\ast$.  We summarize the change point locations in a vector
\begin{align*}
J(\vartheta^{\ast})=(\tau_1^{\ast},\ldots,\tau_{K^{\ast}}^{\ast})
\end{align*}
of dimension   $\left|J(\vartheta^{\ast})\right|$.  
For the sake  of simplicity  we restrict ourselves to estimators  
of the form  $\hat \vartheta  (\cdot ) =\sum_{k=0}^{\hat{K}} \hat \theta_k \mathbb{1}_{[\hat\tau_k,\hat \tau_{k+1})} (\cdot )$
where the estimates $\hat\tau_k$  of the change point locations only attain values 
 at the sampling points $0, \frac 1n, \ldots \frac{n-1}n , 1$ and denote the  set of these functions  by  $\mathcal{S}_n$. 
 Following  \cite{frick2014} we propose  to  test 
 for a candidate step function $\vartheta  (\cdot ) =\sum_{k=0}^{K} \theta_k \mathbb{1}_{[\tau_k,\tau_{k+1})} (\cdot )  \in \mathcal{S}_n $
on each interval $[i/n,j/n]$ where $\vartheta$ is constant whether  $\vartheta^\ast$ is constant on this interval as well with  the same  value as $\vartheta$.
For this purpose we use the multiscale statistic
\begin{align}  \label{vn}
V_n(Y,\vartheta)=\max_{0\leq k\leq K} \max_{\substack{n\tau_k\leq i\leq j<n\tau_{k+1} \\ j-i+1\geq nc_n}} \left\{ \frac{1}{\hat \sigma_\star}\sqrt{j-i+1} \left|\overline{Y}_i^j-\theta_k\right|-\sqrt{2\log\frac{en}{j-i+1}}\right\},
\end{align}
 where $\{c_n\}_{n\in \N}$ is a  positive  sequence converging to $0$,
\begin{align*}
\overline{Y}_i^j:=\frac{1}{j-i+1}\sum_{\ell=i}^j Y_{\ell}
\end{align*}
is a local mean  and 
$\hat \sigma_\star^{2}$ is an appropriate estimator of the long run variance  \eqref{lrv}, which will be defined later.
The estimator of the piecewise constant  function $\vartheta^\ast$
is then required to minimize the number of change points over the acceptance region of this multiscale test.  More precisely, 
for a fixed   threshold $ q$  chosen  according to the (asymptotic) null distribution of $V_n$  the  step function estimator $\hat{\vartheta}$ 
is required to fulfil a data fit claim of the form 
\begin{align*}
V_n(Y,\hat{\vartheta})\leq q~,
\end{align*}
and to satisfy simultaneously a parsimony requirement concerning its number of change points. 
This is achieved by first  estimating the number of change points  $K^\ast$ by
\begin{align*}
\hat{K}=\hat{K}(V_n,q)=\inf_{\substack{\vartheta\in\mathcal{S}_n\\ V_n(Y,\vartheta)\leq q}} |J(\vartheta)|.
\end{align*}
Next, we identify among all  suitable candidate step functions  the one which provides the best fit to the data, that is 
\begin{align} \label{schaetzer_4.1}
\hat{\vartheta}=\argmin_{\vartheta\in\mathcal{C}(V_n,q)}\sum_{i=1}^n\Big (Y_i-\vartheta\big (\tfrac{i}{n}\big )\Big )^2~,
\end{align}
where 
\begin{align*}
\mathcal{C}(V_n,q):=\{\vartheta\in\mathcal{S}_n \,:\,|J(\vartheta)|=\hat{K} \text{ and } V_n(Y,\vartheta)\leq q\}
\end{align*}
is  a ``confidence set'' of all functions in $ \mathcal{S}_n$ satsifying the multiscale criterion with a minimal number of change points.
The estimator can  be efficiently computed by a dynamic program and is implemented 
with the function \textit{stepFit} in the R-package \textit{stepR} [see \cite{rpackage}].

The appropriate estimation of the long run variance  $\sigma_\star^2$  is crucial for a good  performance of SMUCE  if it is applied to correlated data, and for this purpose 
we propose a two step procedure as considered in \cite{wuzhao2007}.  We divide the sample in  $m_n=\lfloor \frac{n}{k_n}\rfloor$ blocks
$\{ Y_{1}, \ldots , Y_{k_{n}}\}$, $\{ Y_{k_{n}+1}, \ldots , Y_{2k_{n}}\} , \ldots  , \{ Y_{(m_{n}-1)k_{n}+1}, \ldots , Y_{m_{n}k_{n}}\}  $ of length $k_{n}$ and 
calculate  local averages
\begin{align*} 
A_i:=\frac{1}{k_n}\sum_{j=1}^{k_n} Y_{j+ik_n},
\end{align*}
to mimic the dependence structure of the data. Secondly, we use the difference based estimate
\begin{align}  \label{lrvest}
\hat{\sigma}_\star^2:=\frac{k_n}{2(m_n-1)}\sum_{i=1}^{m_n-1} \left|A_i-A_{i-1}\right|^2,
\end{align}
to eliminate the signal.  Here $k_n$  increases with the sample size in order to achieve the correct asymptotic behaviour. For details see Proposition \ref{estimate_rate} below, where
we prove the consistency of this estimate. 

 In the Gaussian case the only difference to the SMUCE procedure regards the use of the long run variance estimator. Note, however, that we will discuss
  arbitrary dependent error processes, not necessarily Gaussian, in which case the asymptotic analysis of the procedure is substantially more difficult.
   This analysis will be carefully carried out  in the following Section \ref{sec3}. The finite sample properties of the
new multiscale method are investigated by means of a simulation study in Section \ref{sec4}.

\section{Asymptotic properties}
 \label{sec3}
\def\theequation{3.\arabic{equation}}
\setcounter{equation}{0}

Consider the location scale model  \eqref{model_1} with a stationary  error process $\varepsilon= \{\varepsilon_i\}_{i\in\Z}$  such that 
$\E\left[\varepsilon_i\right]=0$, $\Var\left[\varepsilon_i\right]=\sigma^2 > 0 $. For the asymptotic analysis of the multiscale procedure introduced in Section \ref{sec2}
we  assume that $\varepsilon$ is a physical system as introduced in \cite{wu2005}. This means that there exists
a sequence of independent identically distributed  random variables $\{\eta_i\}_{i\in\Z}$ with values in some measure space $\mathcal{S}$ and a measurable function $G:\mathcal{S}^{\mathbb{N}}\rightarrow\R$
such that for all $i  \in \Z$
\begin{align*}
\varepsilon_i=G(\ldots,\eta_{i-1},\eta_i) ~.
\end{align*}
As pointed out  by \cite{wu2011}, physical systems include many of the commonly used  time series models such as  ARMA and GARCH processes.

In the following discussion let $p \geq  1$   and  define  for   a random variable $X$ (in the case of its existence)  $|| X ||_p =  \big (  \mathbb{E} [ |X|^{p } ]\big )^{1/p} $.   
If $\|\varepsilon_i  \| _{p}<  \infty  $ we consider   the physical dependence  measure
\begin{align*}
\delta_{i,p}:=||\varepsilon_i-\varepsilon_i^\star||_p,
\end{align*}
where the random variable $\varepsilon_i^\star$ is defined by   $\varepsilon_i^\star=G(\ldots,\eta_{-1},\eta_0^\prime,\eta_1,\ldots,\eta_i)$  and   $\eta_0^\prime$ is an independent copy of $\eta_0$. 
We also define the quantity 
\begin{align*}
\Delta_{m,p}:=\sum_{i=m}^\infty \delta_{i,p}, \quad m=1,2,\ldots
\end{align*} 
and call a system  $ \{\varepsilon_i\}_{i\in\Z}$  $p$-strong stable if  $\Delta_{0,p}<\infty$  [see \cite{wu2005}]. 
It can be shown that for a  $2$-strong stable   process $\{\varepsilon_i\}_{i\in\Z}$  the covariance function is absolutely summable and thus the long run variance in \eqref{lrv} exists [see e.g.\ \cite{wuphoumaradi2009}].
A further quantity that we will make use of is the so-called projection operator, which for $i\in\Z$ is given by
\begin{align*}
P_i\,\cdot :=\E\left[\cdot\,|\,\mathcal{F}_i\right]-\E\left[\cdot\,|\,\mathcal{F}_{i-1}\right],
\end{align*}
where $\mathcal{F}_i=(\ldots,\eta_{i-1},\eta_i)$. It is shown  in  \cite{wu2011}
that for a  $2$-strong stable   process $\{\varepsilon_i\}_{i\in\Z}$
 the long run variance   \eqref{lrv}
can be represented as $\sigma_\star^2=\E[(\sum_{j=0}^\infty P_0 \varepsilon_j)^2]$.

For the statement of the asymptotic properties in this section we will make the following basic assumptions
\begin{enumerate}
\item[(A1)] $\| \varepsilon_i\|_{4} < \infty $
\item[(A2)] $\Delta_{0,4}<\infty$  and  $\sum_{i=1}^{\infty} i\delta_{i,2}<\infty$
\item[(A3)] $\Delta_{m,3}=\mathcal{O}(m^{-\gamma})$ for some $\gamma>0$ 
\end{enumerate}
Assumption (A3) is used to construct a    simultaneous Gaussian approximation  of the partial sums of the errors
$\varepsilon_i$
 (see Section \ref{sec5} for details).
Assumption (A2) is needed for a proof of the first result of this section, which establishes the 
consistency  of the estimator \eqref{lrvest}  for  the long run variance  with an explicit rate.  
For its precise statement we introduce  the notation  $a_n\asymp b_{n }$ for two sequences $\{a_{n}\}_{n\in \N} $ and $\{b_{n}\}_{n\in \N} $, which 
 means that
\begin{align*}
0<\liminf_{n\rightarrow\infty} \left|a_{n}/ b_n\right|\leq \limsup_{n\rightarrow\infty} \left| a_{n}/ b_n \right|<\infty.
\end{align*}

\begin{prop} \label{estimate_rate}
Consider the nonparametric regression model \eqref{model_1} 
with a piecewise constant regression function \eqref{21}. If   assumptions (A1) and (A2) are satisfied 
and  $k_n\asymp n^{1/3}$, we have for the estimator in \eqref{lrvest}
\begin{align*}
\hat{\sigma}_\star-\sigma_\star=\mathcal{O}_{\mathbb{P}}\left(n^{-1/3} \right),
\end{align*}
where $\sigma_\star^2$ is the long run variance in \eqref{lrv}. 
\end{prop}

Throughout this paper we will always assume that $k_n\asymp n^{1/3}$, if the long run variance estimator  \eqref{lrvest} is used.
Our first  main result shows that the asymptotic null distribution of  the statistic $V_{n}$ 
does not change in the case of dependent observations.

\begin{theorem} \label{Thm1}
Consider the nonparametric regression model \eqref{model_1} 
 with  piecewise constant regression function \eqref{21}.
If  assumptions (A1)--(A3) are satisfied with $\gamma>1/2$ in (A3), $c_{n}\rightarrow 0$ and 
\begin{align} \label{schneller_als}
\lim_{n\rightarrow\infty} \frac{\left(\log n\right)^3}{n^{m(\gamma)} c_n}=0 ~,
\end{align} 
where $ m(\gamma)=\frac{2\gamma-1}{1+6\gamma}$, then it holds 
\begin{align*} 
V_{n}(Y,\vartheta^\ast)\dn \max_{0\leq k\leq K^\ast} \sup_{\tau_k^\ast\leq s < t\leq \tau_{k+1}^\ast} \left \{\frac{|B(t)-B(s)|}{\sqrt{t-s}}-\sqrt{2 \log \frac{e}{t-s}}\right \} \text{ as } n\rightarrow\infty,
\end{align*}
where $\{B(t)\}_{t\in[0,1]}$ denotes a standard Brownian motion.
\end{theorem}

With exactly the same arguments as given in \cite{frick2014}, we can assure for given $\alpha\in(0,1)$ that
\begin{align} \label{alpha}
\lim_{n\rightarrow\infty} \mathbb{P}\left(\hat{K}(V_{n},q)>K^\ast\right)\leq \alpha,
\end{align}
where $q$ is chosen as the $(1-\alpha)$--quantile of 
\begin{align} \label{M}
M:=\sup_{0\leq s\leq t\leq 1} \left \{\frac{|B(t)-B(s)|}{\sqrt{t-s}}-\sqrt{2\log\frac{e}{t-s}}\right\}.
\end{align}
Note that  the distribution of $M$ coincides with the asymptotic distribution in Theorem \ref{Thm1}, if the function $\vartheta^\ast$ is constant (that  is $K^\ast =0$).
We also obtain from Theorem \ref{Thm1} and the  definition of $\hat{K}$ 
 that the probability of overestimating the number of change points becomes arbitrarily small with an increasing sample size.

\begin{corollary} \label{overestimating}
If  the assumptions from Theorem \ref{Thm1}  are satisfied and $q_n\rightarrow\infty$, we have
\begin{align*}
\lim_{n\rightarrow\infty} \mathbb{P}\left(\hat{K}(V_{n},q_n)>K^\ast\right)=0.
\end{align*}
\end{corollary}

The following result shows that the probability of underestimating the true number of change points also converges to $0$  for an increasing sample size.

\begin{theorem} \label{Thm2}
If  the assumptions from Theorem \ref{Thm1} hold and the sequence $\{q_n\}_{n\in\N}$ fulfils
\begin{align} \label{assumption_concerning_q_n}
q_n=\emph{o}\left(\sqrt{n}\right)  
\end{align}
as $n\rightarrow\infty$,  then it follows that
\begin{align*}
\lim_{n\rightarrow \infty} \mathbb{P}\left(\hat{K}(V_{n},q_n)<K^\ast\right)= 0.
\end{align*}
\end{theorem}

Combining Corollary \ref{overestimating} and Theorem \ref{Thm2} yields model selection consistency. 

\begin{corollary} 
If the same assumptions as in Theorem \ref{Thm2} are satisfied  and   $q_n\rightarrow\infty$, then it follows that
\begin{align*}
\lim_{n\rightarrow\infty} \mathbb{P}\left(\hat{K}\left(V_{n},q_n\right)=K^\ast\right)= 1.
\end{align*}
\end{corollary}

Under appropriate assumptions, the change point locations of $\vartheta^\ast$ are estimated correctly.
More precisely, we have the following result.

\begin{theorem} \label{Thm3}
If the assumptions from Theorem \ref{Thm1} hold and the sequence $\{q_n\}_{n\in\N}$ additionally fulfils $q_n\rightarrow \infty$ and 
\begin{align} \label{assumptions_concerning_qn_2}
q_n+\sqrt{2\log\frac{e}{c_n}}=\emph{o}\left(\sqrt{nc_n}\right)  
\end{align}
as $n\rightarrow\infty$, it follows that
\begin{align*}
\lim_{n\rightarrow\infty}\mathbb{P}\Big (\sup_{\vartheta\in\mathcal{C}(V_{n},q_n)} \max_{\tau^\ast\in J(\vartheta^\ast)} \min_{\tau\in J(\vartheta)} \left|\tau^\ast-\tau\right|> c_n\Big)=0.
\end{align*}
In particular we have for $k=1,\ldots,K^\ast$
\begin{align*}
\lim_{n\rightarrow\infty} \mathbb{P}\Big (\sup_{\vartheta\in\mathcal{C}(V_{n},q_n)} \left|\tau_k^\ast-\tau_k\right|> c_n\Big)= 0.
\end{align*}
\end{theorem}

\section{Finite sample properties}
 \label{sec4}
\def\theequation{4.\arabic{equation}}
\setcounter{equation}{0}
In this section we  compare the finite sample performance of the  change point estimator developed and analyzed in Section \ref{sec3} with SMUCE and the change point estimator proposed by  \cite{tecuapetla2017}  for $m$-dependent errors.
 These authors use the abbreviation  JUSD for their procedure and we will use 
the notation   DepSMUCE for the procedure  \eqref{schaetzer_4.1} developed in this paper. The sample size is $n=1000$ and  all  results are based on  $1000$ simulation runs. For DepSMUCE, we consider a block length of $k=10$. 
Concerning the change point estimator JUSD, it is necessary to specify a value for $m$. The R-package \textit{dbacf} [see \cite{rpackage2}] provides a graphical procedure to choose $m$
which is used throughout the simulation study. 

We compare the deviations between the estimated and the true number of change points, and the mean deviation of $|K^\ast-\hat{K}|$. Concerning the data fit, we compute the  mean squared error
 
\begin{align*}
\text{MSE}(\hat{\vartheta}):=\frac{1}{n}\sum_{i=1}^n \left(\vartheta^\ast\left(\tfrac{i}{n}\right)-\hat{\vartheta}\left(\tfrac{i}{n}\right)\right)^2
\end{align*}
and mean absolute deviation 
\begin{align*}
\text{MAE}(\hat{\vartheta}):=\frac{1}{n}\sum_{i=1}^n \left|\vartheta^\ast\left(\tfrac{i}{n}\right)-\hat{\vartheta}\left(\tfrac{i}{n}\right)\right|,
\end{align*}
respectively.
Furthermore, we  also present  histograms of the estimated locations of the changes for all three estimators. 

All procedures depend sensitively on the  threshold $q$ in the definition of the change point estimator and we investigate three different choices of $q$. More precisely, considering \eqref{alpha}, we choose the significance level $\alpha$ as $0.1$, $0.5$ and $0.9$ and set $q$ as the $(1-\alpha)$--quantile of the distribution of the random variable $M$ in \eqref{M}. Since this quantile cannot be derived directly, we perform Monte Carlo simulations of the test statistic $V_{n}(Y,\vartheta^\ast)$ with $\vartheta^\ast\equiv 0$ and independent standard normal distributed errors, i.e. $\varepsilon_i\sim\mathcal{N}(0,1)$ (based on $10000$ repetitions). This is exactly the same procedure as in the R-package \textit{stepR} [see \cite{rpackage}]. 

First, we illustrate that SMUCE is relatively robust to weak dependencies but it does not yield satisfactory results when the innovations exhibit a stronger dependence. To this end, we consider two MA($1$) error processes with different MA parameters. Let
\begin{align} \label{MA1_process}
\varepsilon_i=\eta_i+\kappa \eta_{i-1}, \quad i\in\Z,
\end{align}
where  $\{\eta_i\}_{i\in \Z}$  is a sequence of standard normal distributed errors. We consider the cases $\kappa=0.1$ and $\kappa=0.3$, respectively, and assume that the function $\vartheta^\ast$ in model \eqref{model_1}  
has  $K^\ast=5$ change points at  locations 
\begin{align} \label{change_point_locations}
(\tau_1^\ast,\tau_2^\ast,\tau_3^\ast,\tau_4^\ast,\tau_5^\ast)=(101/1000,301/1000,501/1000,551/1000,751/1000).
\end{align}
The corresponding function intensities are given by
\begin{align} \label{function_values_MA1}
(\theta_0^\ast,\theta_1^\ast,\theta_2^\ast,\theta_3^\ast,\theta_4^\ast,\theta_5^\ast)=(0,1,0,2,0,-1).
\end{align}

\begin{table}[!htbp]
\centering
{\footnotesize
\begin{tabular}{l|l|l|l|l|l|l|l} 
$\hat{K}-K^\ast$ & $\leq -3$ & $-2$ & $-1$ & $0$ & $+1$ & $+2$ & $\geq +3$ \\
\hline
 SMUCE(0.1) & 0.000 & 0.000 & 0.001 & 0.980 & 0.019 & 0.000 & 0.000\\
 SMUCE(0.5) & 0.000 & 0.000 & 0.000 & 0.760 & 0.209 & 0.031 & 0.000\\
 SMUCE(0.9) & 0.000 & 0.000 & 0.000 & 0.238 & 0.343 & 0.267 & 0.152\\
\hline
 DepSMUCE(0.1) & 0.000 & 0.000 & 0.117 & 0.883 & 0.000 & 0.000 & 0.000 \\
 DepSMUCE(0.5) & 0.000 & 0.000 & 0.009 & 0.988 & 0.003 & 0.000 & 0.000\\
 DepSMUCE(0.9) & 0.000 & 0.000 & 0.000 & 0.946 & 0.053 & 0.001 & 0.000 \\
\hline
 JUSD(0.1) & 0.075 & 0.065 & 0.125 & 0.702 & 0.027 & 0.004 & 0.003 \\
 JUSD(0.5) & 0.020 & 0.024 & 0.080 & 0.745 & 0.087 & 0.020 & 0.023 \\
 JUSD(0.9) & 0.003 & 0.004 & 0.033 & 0.631 & 0.168 & 0.085 & 0.076 \\
 \hline
\end{tabular}
\caption{\it Proportion  of estimated numbers of change points  (the true number of change points is $K^\ast=5$)  in model \eqref{model_1}
with step function defined by \eqref{change_point_locations} and \eqref{function_values_MA1} and  an MA($1$) error process defined  in \eqref{MA1_process} with $\kappa=0.1$.}
\label{MA1_01}
}
\end{table}

\begin{table}[!htbp]
\centering
{\footnotesize
\begin{tabular}{l|l|l|l|l|l|l|l} 
$\hat{K}-K^\ast$ & $\leq -3$ & $-2$ & $-1$ & $0$ & $+1$ & $+2$ & $\geq +3$ \\
\hline
 SMUCE(0.1) & 0.000 & 0.000 & 0.001 & 0.619 & 0.302 & 0.066 & 0.012\\
 SMUCE(0.5) & 0.000 & 0.000 & 0.000 & 0.069 & 0.184 & 0.262 & 0.486\\
 SMUCE(0.9) & 0.000 & 0.000 & 0.000 & 0.000 & 0.010 & 0.025 & 0.965\\
\hline
 DepSMUCE(0.1) & 0.001 & 0.043 & 0.356 & 0.600 & 0.000 & 0.000 & 0.000 \\
 DepSMUCE(0.5) & 0.000 & 0.000 & 0.048 & 0.947 & 0.005 & 0.000 & 0.000\\
 DepSMUCE(0.9) & 0.000 & 0.000 & 0.007 & 0.919 & 0.070 & 0.004 & 0.000 \\
\hline
 JUSD(0.1) & 0.231 & 0.124 & 0.179 & 0.446 & 0.018 & 0.001 & 0.001 \\
 JUSD(0.5) & 0.072 & 0.077 & 0.167 & 0.618 & 0.043 & 0.016 & 0.007 \\
 JUSD(0.9) & 0.010 & 0.016 & 0.111 & 0.615 & 0.156 & 0.060 & 0.032 \\
 \hline
\end{tabular}
\caption{\it Proportion  of estimated numbers of change points  (the true number of change points is $K^\ast=5$)  in model \eqref{model_1}
with  step function defined by \eqref{change_point_locations} and \eqref{function_values_MA1}  and  an MA($1$) error process defined  in \eqref{MA1_process} with $\kappa=0.3$.}
\label{MA1_03}
}
\end{table}

Tables \ref{MA1_01} and \ref{MA1_03} show the performance of SMUCE, DepSMUCE and JUSD in the estimation of the number of change points for different values of $\alpha$. For example,
in Table \ref{MA1_01} we display results for model \eqref{MA1_process} with  $\kappa=0.1$ and we observe that DepSMUCE estimates the correct number of change points in $98.8\%$ of the cases
if we work with  $\alpha=0.5$. Considering the first three rows of Table \ref{MA1_01}, it can be seen that SMUCE performs relatively well if $\alpha=0.1$. However for $\alpha=0.5$ DepSMUCE  already shows some improvement and 
for  $\alpha=0.9$ DepSMUCE and JUSD  show a better performance because they are constructed to  address  dependency in the data.
The advantages of  these  two procedures become even more visible  in Table \ref{MA1_03}, where we consider  a stronger dependence, that is 
$\kappa=0.3$. In this case SMUCE tends to overestimate the true number of change points.
We also observe  a better performance of DepSMUCE  compared to JUSD, which often estimates a too small number of change points.
\begin{figure} [ht]
  \centering
{\includegraphics[width=0.3\textwidth, height=160px]{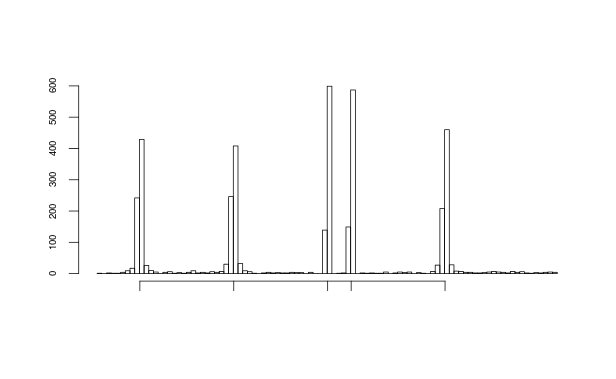}}\qquad
{\includegraphics[width=0.3\textwidth, height=160px]{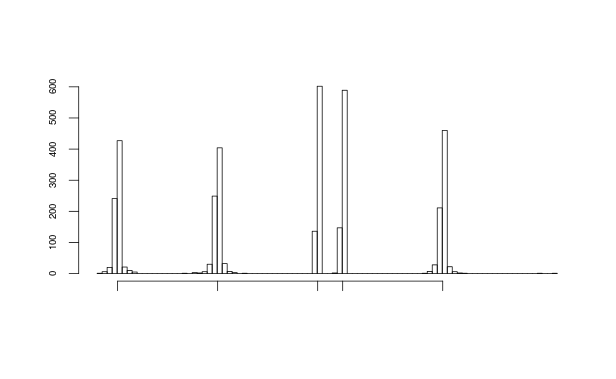}}\qquad
{\includegraphics[width=0.3\textwidth, height=160px]{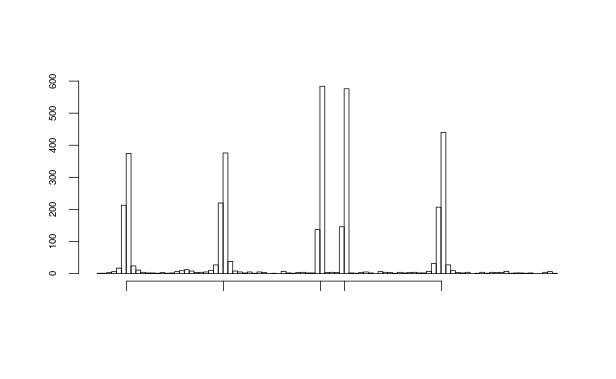}} \\
\par
\vspace{-1.2cm}
{\includegraphics[width=0.3\textwidth, height=160px]{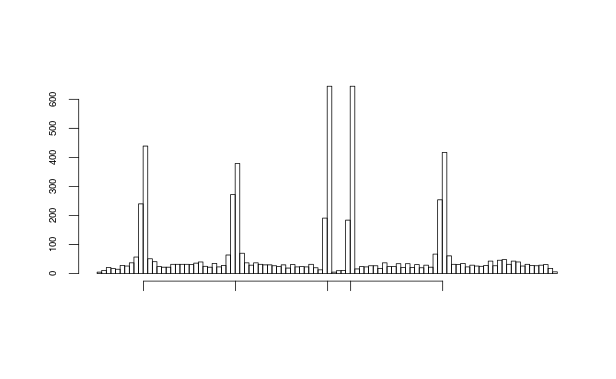}}\qquad
{\includegraphics[width=0.3\textwidth, height=160px]{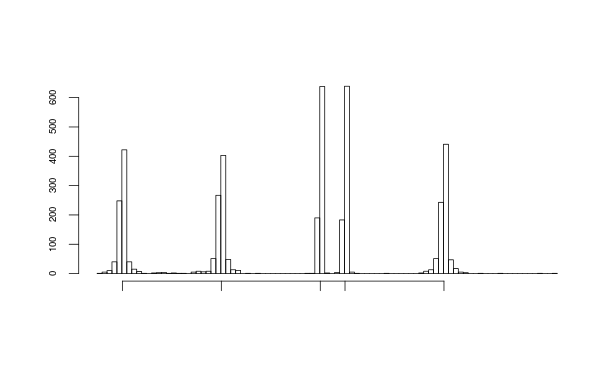}}\qquad
{\includegraphics[width=0.3\textwidth, height=160px]{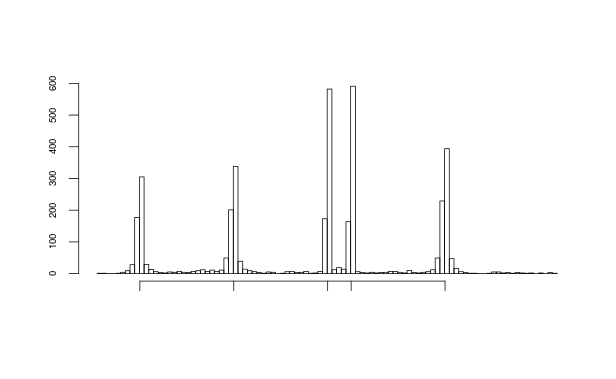}}
\vspace{-1.11cm}
  \caption{\it  
  Histograms of estimated change point locations for different estimators. First row: MA$(1)$  error process with $\kappa=0.1$. Second row: MA$(1)$ error process with $\kappa=0.3$. Left 
  column:  SMUCE.
 Middle  column: DepSMUCE.  Right  column: JUSD. The ``true''  change points are located at  $101$, $301$, $501$, $551$, and $751$.}
   \label{HistogramsMA1}
\end{figure}
Our findings are confirmed by Table \ref{MA1_01_MSE}, where we display the average MSE, MAE, and the average value of $|K^\ast-\hat{K}|$. Figure \ref{HistogramsMA1} shows histograms of the estimated change point locations for $\alpha=0.5$. The comparatively bad performance of JUSD can be explained by the fact that it requires the specification of the order of the MA($m$) process. We observed that the data driven procedure to choose $m$ from the R-package \textit{dbacf} [see \cite{rpackage2}] does not work well for small MA parameters. A simulation study with $\kappa=0.5$, which is not included here for the sake of brevity, shows that JUSD works much better for 
larger MA parameters, and as a consequence  the behaviour of DepSMUCE and JUSD becomes more similar.

\begin{table}[ht]
\centering
{\footnotesize
\begin{tabular}{l||l|l|l||l|l|l}
& \multicolumn{3}{l||}{\phantom{=====}$\kappa=0.1$} & \multicolumn{3}{l}{\phantom{=====}$\kappa=0.3$}\\
 \cline{2-7}
 & $|K^\ast-\hat{K}|$ & MSE & MAE & $|K^\ast-\hat{K}|$ & MSE & MAE\\
\hline
SMUCE(0.1) & 0.020 & 0.018 & 0.060 & 0.475 & 0.033 & 0.093\\
SMUCE(0.5) & 0.271 & 0.019 & 0.064 & 2.569 & 0.045 & 0.118\\
SMUCE(0.9) & 1.407 & 0.024 & 0.077 & 6.488 & 0.063 & 0.145\\
\hline
DepSMUCE(0.1) & 0.117 & 0.025 & 0.072 & 0.446 & 0.064 & 0.139 \\
DepSMUCE(0.5) & 0.012 & 0.018 & 0.060 & 0.053 & 0.031 & 0.088\\
DepSMUCE(0.9) & 0.056 & 0.018 & 0.060 & 0.084 & 0.030 & 0.085\\
\hline
JUSD(0.1) & 0.549 & 0.050 & 0.109 & 1.226 & 0.117 & 0.209\\
JUSD(0.5) & 0.397 & 0.031 & 0.082 & 0.647 & 0.069 & 0.145\\
JUSD(0.9) & 0.705 & 0.024 & 0.073 & 0.577 & 0.044 & 0.109 \\
 \hline
\end{tabular}
\caption{\it
Average of $|K^\ast-\hat{K}|$,  MSE, and  MAE
of different estimates in  the MA$(1)$-model \eqref{MA1_process}.}
\label{MA1_01_MSE}
}
\end{table}

This observation is also confirmed in our next example, where we consider an MA($4$) error process with relatively large parameters, that is
\begin{align} \label{MA_4_process}
\varepsilon_i=\eta_i+0.9\eta_{i-1}+0.8\eta_{i-2}+0.7\eta_{i-3}+0.6\eta_{i-4}, \quad i\in\Z ~.
\end{align}
Here $\{\eta_i\}_{i \in \Z} $ denotes again a sequence of independent standard normal distributed errors. We assume that the function $\vartheta^\ast$ in model \eqref{model_1}  
has  $K^\ast=5$ change points at the locations given in \eqref{change_point_locations} and that the corresponding function intensities are given by 
\begin{align} \label{function_values_MA4}
(\theta_0^\ast,\theta_1^\ast,\theta_2^\ast,\theta_3^\ast,\theta_4^\ast,\theta_5^\ast)=(0,3,0,4,0,-3).
\end{align}
\begin{table}[ht]
\centering
{\footnotesize
\begin{tabular}{l|l|l|l|l|l|l|l} 
$\hat{K}-K^\ast$ & $\leq -3$ & $-2$ & $-1$ & $0$ & $+1$ & $+2$ & $\geq +3$ \\
\hline
 SMUCE(0.1) & 0.000 & 0.000 & 0.000 & 0.000 & 0.000 & 0.000 & 1.000\\
 SMUCE(0.5) & 0.000 & 0.000 & 0.000 & 0.000 & 0.000 & 0.000 & 1.000\\
 SMUCE(0.9) & 0.000 & 0.000 & 0.000 & 0.000 & 0.000 & 0.000 & 1.000\\
\hline
 DepSMUCE(0.1) & 0.020 & 0.138 & 0.511 & 0.330 & 0.001 & 0.000 & 0.000 \\
 DepSMUCE(0.5) & 0.000 & 0.006 & 0.173 & 0.806 & 0.015 & 0.000 & 0.000\\
 DepSMUCE(0.9) & 0.000 & 0.000 & 0.024 & 0.856 & 0.114 & 0.006 & 0.000 \\
\hline
 JUSD(0.1) & 0.025 & 0.188 & 0.511 & 0.275 & 0.001 & 0.000 & 0.000 \\
 JUSD(0.5) & 0.000 & 0.006 & 0.145 & 0.812 & 0.037 & 0.000 & 0.000 \\
 JUSD(0.9) & 0.000 & 0.000 & 0.013 & 0.709 & 0.240 & 0.032 & 0.006 \\
 \hline
\end{tabular}
}
\caption{\it Proportion  of estimated numbers of change points  (the true number of change points is $K^\ast=5$)  in model \eqref{model_1}
with  step function defined by \eqref{change_point_locations} and \eqref{function_values_MA4}  and  an MA($4$) error process defined  in \eqref{MA_4_process}.}
\label{MA_4}
\end{table}

\begin{table}[ht]
\centering
{\footnotesize
\begin{tabular}{l||l|l|l||l|l|l}
& \multicolumn{3}{l||}{\phantom{===i==}MA(4)} & \multicolumn{3}{l}{\phantom{====}ARMA(2,6)}\\
\cline{2-7}
 & $|K^\ast-\hat{K}|$ & MSE & MAE & $|K^\ast-\hat{K}|$ & MSE & MAE\\
\hline
SMUCE(0.1)& 43.845 & 1.787 & 1.016 & 59.950 & 4.174 & 1.592\\
SMUCE(0.5) & 56.842 & 2.041 & 1.108 & 73.800 & 4.550 & 1.684\\
SMUCE(0.9) & 67.865 & 2.208 & 1.166 & 85.582 & 4.798 & 1.743\\
\hline
DepSMUCE(0.1) & 0.848 & 0.861 & 0.584 & 0.516 & 1.534 & 0.778  \\
DepSMUCE(0.5) & 0.200 & 0.418 & 0.364 & 0.064 & 0.646 & 0.465\\
DepSMUCE(0.9) & 0.150 & 0.319 & 0.322 & 0.115 & 0.586 & 0.449\\
\hline
JUSD(0.1) & 0.963 & 0.929 & 0.619  & 1.422 & 1.720 & 0.879\\
JUSD(0.5) & 0.194 & 0.401 & 0.357 & 2.181 & 1.042 & 0.637 \\
JUSD(0.9) & 0.336 & 0.343 & 0.341& 3.979 & 1.019 & 0.630  \\
 \hline
\end{tabular}
\caption{\it Average of $|K^\ast-\hat{K}|$,  MSE, and  MAE
of different estimates
 under the same model assumptions as in Table \ref{MA_4} and Table \ref{ARMA2,6}.}
\label{MA_4_MSE,ARMA26_MSE}
}
\end{table}

\begin{figure} [ht]
  \centering
    {\includegraphics[width=0.3\textwidth, height=160px]{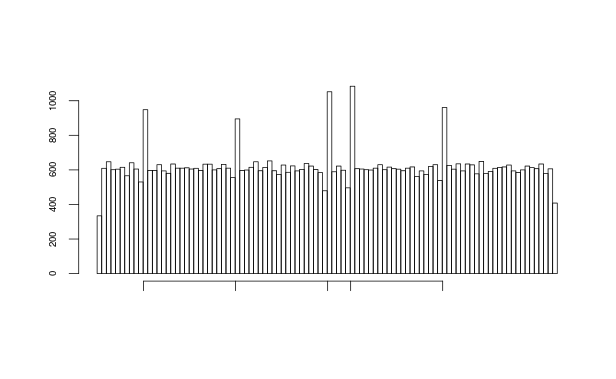}}\qquad
   {\includegraphics[width=0.3\textwidth, height=160px]{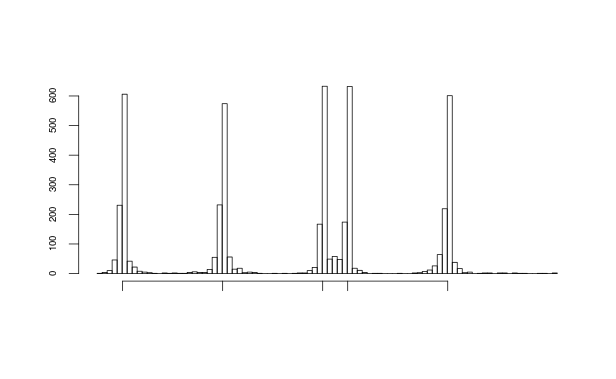}}\qquad
{\includegraphics[width=0.3\textwidth, height=160px]{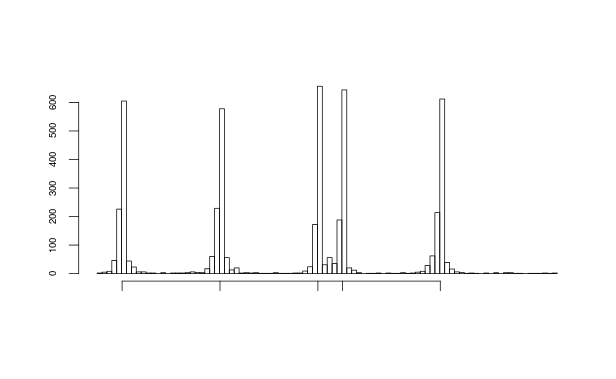}}\\
\par
\vspace{-1.2cm}
{\includegraphics[width=0.3\textwidth, height=160px]{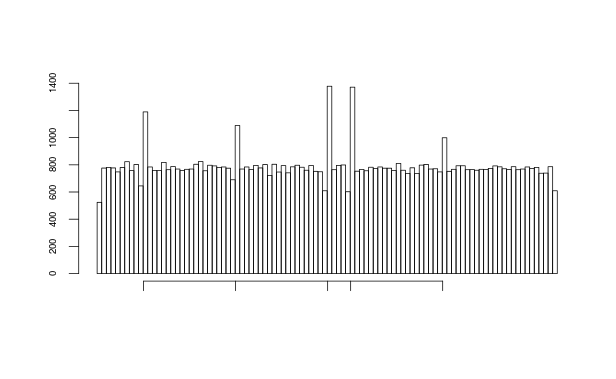}}\qquad
{\includegraphics[width=0.3\textwidth, height=160px]{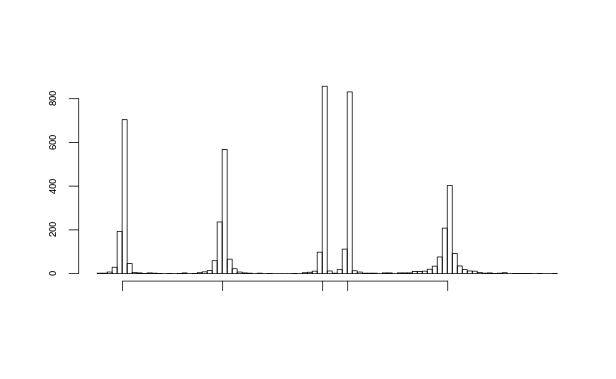}}\qquad
{\includegraphics[width=0.3\textwidth, height=160px]{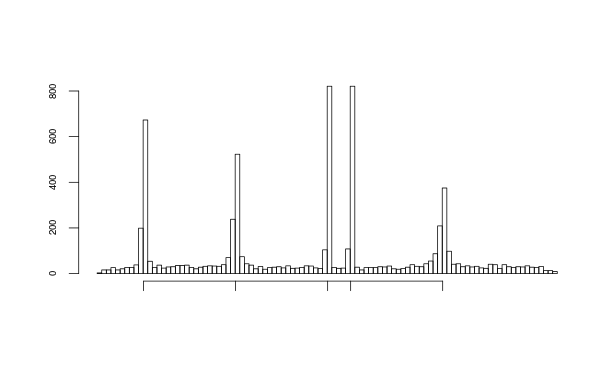}}
\vspace{-1.11cm}
  \caption{\it  
  Histograms of estimated change point locations for different estimators.
  Upper row: MA$(4)$ error process. Lower row: ARMA$(2,6)$  error process. Left columun:  SMUCE.
 Middle columun: DepSMUCE.  Right columun: JUSD. The ``true''  change points are located at  $101$, $301$, $501$, $551$, and $751$.}
   \label{Histograms2}
\end{figure}

The data driven rule in the R-package \textit{dbacf}  works well and determines $m=4$ for JUSD correctly in all of the iterations. 
Table \ref{MA_4} shows the performance  
of  SMUCE, DepSMUCE, and JUSD. For example, if $\alpha=0.5$, DepSMUCE estimates the number $K^\ast$ of change points correctly in 
$80.6\%$ of the cases, while it underestimates $K^\ast$ by $1$ in $17.3\%$ of the cases.
The first three rows show that SMUCE is not able to correctly estimate the
 number of change points in the case of an MA($4$) error process. Of course, SMUCE  is designed for independent data, but  it always estimates a much larger number of change points than $5$.
  In contrast, JUSD and DepSMUCE perform substantially better if they are used with $\alpha=0.5$.
 In particular, they yield very similar results and DepSMUCE is able to compete with JUSD, which is specially designed for  $m$-dependent processes
 (note that  we used the correct $m$ in the simulations). Similar observations can be made for the estimation error  (see the left part of Table \ref{MA_4_MSE,ARMA26_MSE}).
  These observations are also supported by the upper part of Figure \ref{Histograms2} which   shows the histograms of the estimated change point locations.
 DepSMUCE and JUSD are able to identify the locations correctly in most of the cases and show a rather similar behaviour.
On the other hand, SMUCE is not reliable for the estimation of the signal in case of strong dependencies. DepSMUCE and  JUSD show a  good and  similar performance if the error process is an MA($4$)-process
and the corresponding parameters are reasonably large.

Finally, we consider an example where the error process  in model \eqref{model_1} is a stationary and causal ARMA($2,6$)-process defined by 
\begin{align} \label{ARMA2,6_process}
\varepsilon_i = 0.75\varepsilon_{i-1}-0.5\varepsilon_{i-2}+\eta_i+0.8\eta_{i-1}+0.7\eta_{i-2}+0.6\eta_{i-3}+
0.5\eta_{i-4}+0.4\eta_{i-5}+0.3\eta_{i-6}, 
\end{align} 
where $\{\eta_i\}_{i \in \Z} $  is a sequence of independent standard normal  distributed random variables.
We consider again a model with  $K^\ast=5$ change points located as described in \eqref{change_point_locations} with the corresponding function intensities
\begin{align} \label{function_values_ARMA2,6}
(\theta_0^\ast,\theta_1^\ast,\theta_2^\ast,\theta_3^\ast,\theta_4^\ast,\theta_5^\ast)=(0,5,1,8,1,-2)
\end{align}
(see Figure \ref{Fig1}). The data driven procedure from the R-package \textit{dbacf} [\cite{rpackage2}] now leads to ambiguous results because there is no correct $m$ to estimate.
Table \ref{ARMA2,6} shows the  estimated numbers of change points. While at  level $\alpha=0.5$ DepSMUCE correctly estimates $K^\ast=5$ in more than $93$\% of the
cases, 
 JUSD  mostly overestimates $K^\ast$. 
 From the right part of Table \ref{MA_4_MSE,ARMA26_MSE}  we also observe  that $K^\ast$ is estimated more precisely by  DepSMUCE
 than by JUSD with smaller MSE and MAE. 
As in the MA($4$)-example, SMUCE in general includes a large amount of false positives. 
Finally, these  results are reflected in Figure \ref{Histograms2}, where we show the histograms of estimated change points.
For the  ARMA($2,6$) error process DepSMUCE yields substantially better  results than JUSD.

\begin{table}[!htbp]
\centering
{\footnotesize
\begin{tabular}{l|l|l|l|l|l|l|l} 
$\hat{K}-K^\ast$ & $\leq -3$ & $-2$ & $-1$ & $0$ & $+1$ & $+2$ & $\geq +3$ \\
\hline
 SMUCE(0.1) & 0.000 & 0.000 & 0.000 & 0.000 & 0.000 & 0.000 & 1.000\\
 SMUCE(0.5) & 0.000 & 0.000 & 0.000 & 0.000 & 0.000 & 0.000 & 1.000\\
 SMUCE(0.9) & 0.000 & 0.000 & 0.000 & 0.000 & 0.000 & 0.000 & 1.000\\
\hline
 DepSMUCE(0.1) & 0.001 & 0.061 & 0.391 & 0.547 & 0.000 & 0.000 & 0.000 \\
 DepSMUCE(0.5) & 0.000 & 0.001 & 0.049 & 0.937 & 0.013 & 0.000 & 0.000\\
 DepSMUCE(0.9) & 0.000 & 0.000 & 0.005 & 0.892 & 0.096 & 0.007 & 0.000 \\
\hline
 JUSD(0.1) & 0.055 & 0.144 & 0.203 & 0.292 & 0.099 & 0.077 & 0.129 \\
 JUSD(0.5) & 0.005 & 0.022 & 0.087 & 0.447 & 0.077 & 0.066 & 0.296 \\
 JUSD(0.9) & 0.001 & 0.002 & 0.012 & 0.444 & 0.067 & 0.033 & 0.441 \\
 \hline
\end{tabular}
\caption{\it  
\it Proportion  of estimated numbers of change points  (the true number of change points is $K^\ast=5$)  in model \eqref{model_1}
with  step function defined by \eqref{change_point_locations} and \eqref{function_values_ARMA2,6}   and  an ARMA $(2,6)$
 error process defined  in \eqref{ARMA2,6_process}.}
\label{ARMA2,6}
}
\end{table}

\section{Proofs}
 \label{sec5}
\def\theequation{5.\arabic{equation}}
\setcounter{equation}{0}

\subsection{Proof of Theorem \ref{Thm1}}

The proof essentially proceeds in two steps.
First we will prove  an analog of the result for the 
statistic 
$$
V_{n,\sigma_\star} (Y,\vartheta)  =\max_{0\leq k\leq K} \max_{\substack{n\tau_k\leq i\leq j<n\tau_{k+1} \\ j-i+1\geq nc_n}} \left\{ \frac{1}{\sigma_\star}\sqrt{j-i+1} \left|\overline{Y}_i^j-\theta_k\right|-\sqrt{2\log\frac{en}{j-i+1}}\right\},
$$
which is defined as $V_n$  using the known long run variance. This result is essentially based on a Gaussian approximation via assumption (A3). In a second step we use Proposition  \ref{estimate_rate}, for which assumptions (A1) and (A2) are needed, to show that the error caused by the estimation  of the long run variance is negligible. The proof of the proposition is given at the end of this section.

\begin{theorem} \label{Thm4}
Consider the nonparametric regression model \eqref{model_1} and assume that $\| \varepsilon_i\|_{3} < \infty$. If assumption (A3) holds for some $\gamma>1/2$ and $c_n\rightarrow 0$ fulfils \eqref{schneller_als}, then we have
\begin{align*} 
V_{n,\sigma_\star} (Y,\vartheta^\ast) \dn  \max_{0 \leq k \leq K^\ast} \sup_{\tau_k^\ast \leq s <t \leq \tau_{k+1}^\ast}  \left \{
\frac {|B(t) - B(s)|}{\sqrt{t-s}} - \sqrt{2 \log \frac {e}{t-s}} \right \} \text{ as } n\rightarrow\infty,
\end{align*}
where $\{B(t)\}_{t\in[0,1]}$ denotes a standard Brownian motion.
\end{theorem}

\underline{\textit{Proof of Theorem \ref{Thm4}} :} 
By (A3), the assumptions of Theorem 1 in \cite{wuzhou2011} are fulfilled. It therefore follows that on a richer probability space $(\check{\Omega},\check{\mathcal{A}},\check{\mathbb{P}})$, there exists a process $\{\check{S}_i\}_{i=1}^n$ and a centered Gaussian process $\{\check{G}_i\}_{i=1}^n$ with independent increments such that
\begin{align} \label{approxi1}
\Big (\sum_{\ell=1}^i \varepsilon_{\ell}\Big )_{i=1}^n \ds \left(\check{S}_i\right)_{i=1}^n \text{ and } \max_{1\leq i\leq n} \left|\check{S}_i-\check{G}_i\right|=\mathcal{O}_{\mathbb{\check{P}}}(\tau_n),
\end{align}
where
\begin{align*}
\tau_n=n^{(1+2\gamma)/(1+6\gamma)}(\log n)^{8\gamma/(1+6\gamma)}.
\end{align*}
Moreover, again on a richer probability space $(\hat{\Omega},\hat{\mathcal{A}},\hat{\mathbb{P}})$, there exists another Gaussian process $\{\hat{G}_i\}_{i=1}^n$ and i.i.d.\ random variables $U_\ell \sim\mathcal{N}(0,\sigma_\star^2)$ such that
\begin{align} \label{approxi2}
\left(\check{G}_i\right)_{i=1}^n\ds \big (\hat{G}_i\big )_{i=1}^n \text{ and } \max_{1\leq i\leq n} \Big |\hat{G}_i-\sum_{\ell=1}^i U_\ell\Big |=\mathcal{O}_{\mathbb{\hat{P}}}(\tau_n).
\end{align}

By \eqref{approxi1}, it follows that $V_{n,\sigma_\star}(Y,\vartheta^\ast)\ds \check{V}_{n,\sigma_\star}(Y,\vartheta^\ast)$, where
\begin{align*}
\check{V}_{n,\sigma_\star}(Y,\vartheta^\ast)=\max_{0\leq k\leq K^\ast} \max_{\substack{n\tau_k^\ast\leq i\leq j<n\tau_{k+1}^\ast \\ j-i+1\geq nc_n}} \left \{\frac{1}{\sigma_\star\sqrt{j-i+1}}\left|\check{S}_j-\check{S}_{i-1}\right|-\sqrt{2\log\frac{en}{j-i+1}}\right \}.
\end{align*}
By \eqref{schneller_als} we have  $\tau_n=\text{o}(\sqrt{nc_n})$, and a further application of  \eqref{approxi1}  and  the triangle inequality 
 yield
\begin{align*}
\check{V}_{n,\sigma_\star}(Y,\vartheta^\ast)=\max_{0\leq k\leq K^\ast} \max_{\substack{n\tau_k^\ast\leq i\leq j<n\tau_{k+1}^\ast \\ j-i+1\geq nc_n}} \left \{\frac{1}{\sigma_\star\sqrt{j-i+1}}\left|\check{G}_j-\check{G}_{i-1}\right|-\sqrt{2\log\frac{en}{j-i+1}}\right \}+\text{o}_{\check{\mathbb{P}}}(1).
\end{align*}

  By \eqref{approxi2}, the first random variable on the right hand side has the same distribution as 
\begin{align*}
\hat{V}_{n,\sigma_\star}(Y,\vartheta^\ast)=\max_{0\leq k\leq K^\ast} \max_{\substack{n\tau_k^\ast\leq i\leq j<n\tau_{k+1}^\ast \\ j-i+1\geq nc_n}} \left \{\frac{1}{\sigma_\star\sqrt{j-i+1}}\left|\hat{G}_j-\hat{G}_{i-1}\right|-\sqrt{2\log\frac{en}{j-i+1}}\right \}.
\end{align*}
With the same arguments as given above, we obtain
\begin{align*}
\hat{V}_{n,\sigma_\star}(Y,\vartheta^\ast)=\max_{0\leq k\leq K^\ast} \max_{\substack{n\tau_k^\ast\leq i\leq j<n\tau_{k+1}^\ast \\ j-i+1\geq nc_n}} \left \{\frac{1}{\sigma_\star}\sqrt{j-i+1}\left|\overline{U}_i^j\right|-\sqrt{2\log\frac{en}{j-i+1}}\right \}+\text{o}_{\hat{\mathbb{P}}}(1).
\end{align*}
Note that
\begin{align*}
&\phantom{=}\max_{0 \leq k \leq K^\ast} \max_{\substack{n\tau_k^\ast\leq i\leq j<n\tau_{k+1}^\ast \\ j-i+1\geq nc_n}} \left \{ \frac{1}{\sigma_\star} \sqrt{j-i+1} \left|\overline{U}_i^j\right| - \sqrt{2\log\frac{en}{j-i+1}} \right \}\\
& \ds \max_{0 \leq k \leq K^\ast} \max_{\substack{n\tau_k^\ast\leq i\leq j<n\tau_{k+1}^\ast \\ j-i+1\geq nc_n}} \left \{\sqrt{j-i+1} \left|\overline{Z}_i^j\right| - \sqrt{2\log\frac{en}{j-i+1}} \right \},
\end{align*}
where $Z_1,\ldots,Z_n\sim \mathcal{N}(0,1)$ are i.i.d.
The assertion now follows with the same arguments as given in the proof of Theorem 1 in \cite{frick2014}.
\hfill $\Box$ \\
\\
\underline{\textit{Proof of Theorem \ref{Thm1}} :}
For the sake of clarity, we will denote  the statistic $V_{n}$ in \eqref{vn} by $V_{n,\hat{\sigma}_{\star}}$ to emphasize its dependence on the estimator  $\hat{\sigma}_\star^2$ of the long run variance. 
Considering the proof of Theorem \ref{Thm4}, it suffices to show that
\begin{align} \label{zuzeigen}
V_{n,\sigma_\star}(Y,\vartheta^\ast)-V_{n,\hat{\sigma}_\star}(Y,\vartheta^\ast)=\text{o}_{\mathbb{P}}(1).
\end{align}
A straightforward application of the  triangle inequality yields
\begin{align*}
\left|V_{n,\sigma_\star}(Y,\vartheta^\ast)-V_{n,\hat{\sigma}_\star}(Y,\vartheta^\ast)\right|\leq \left|\frac{1}{\sigma_\star}-\frac{1}{\hat{\sigma}_\star}\right|\max_{0\leq k\leq K^\ast} \max_{\substack{n\tau_k^\ast\leq i\leq j<n\tau_{k+1}^\ast \\ j-i+1\geq nc_n}} \left| \sqrt{j-i+1}\left(\overline{Y}_i^j-\theta_k^\ast\right)\right|.
\end{align*}
We use again the Gaussian approximation result from Theorem \ref{Thm4} and obtain 
\begin{align*}
D_n:=\max_{\substack{1\leq i\leq j\leq n \\ j-i+1\geq nc_n}} \left|\sqrt{j-i+1}\left(\overline{Y}_i^j-\E[\overline{Y}_i^j]\right)\right|
&\ds \max_{\substack{1\leq i\leq j\leq n \\ j-i+1\geq nc_n}} \Bigl \{\frac{1}{\sqrt{j-i+1}}\left|\check{S}_j-\check{S}_{i-1}\right|\Bigr \}
\end{align*}
as well as
\begin{align*}
\max_{\substack{1\leq i\leq j\leq n \\ j-i+1\geq nc_n}} \Bigl \{\frac{1}{\sqrt{j-i+1}}\left|\check{S}_j-\check{S}_{i-1}\right|\Bigr \}=\max_{\substack{1\leq i\leq j\leq n \\ j-i+1\geq nc_n}} \Bigl \{\frac{1}{\sqrt{j-i+1}}\left|\check{G}_j-\check{G}_{i-1}\right|\Bigr \}+\text{o}_{\check{\mathbb{P}}}(1).
\end{align*}
Furthermore, it holds that
\begin{align*}
\max_{\substack{1\leq i\leq j\leq n \\ j-i+1\geq nc_n}} \Bigl \{\frac{1}{\sqrt{j-i+1}}\left|\check{G}_j-\check{G}_{i-1}\right|\Bigr \}\ds \max_{\substack{1\leq i\leq j\leq n \\ j-i+1\geq nc_n}} \Bigl \{\frac{1}{\sqrt{j-i+1}}\left|\hat{G}_j-\hat{G}_{i-1}\right|\Bigr \}
\end{align*}
and
\begin{align*}
\max_{\substack{1\leq i\leq j\leq n \\ j-i+1\geq nc_n}} \Bigl \{\frac{1}{\sqrt{j-i+1}}\left|\hat{G}_j-\hat{G}_{i-1}\right|\Bigr \}=\max_{\substack{1\leq i\leq j\leq n \\ j-i+1\geq nc_n}} \Bigl\{ \sqrt{j-i+1} \left|\overline{U}_i^j\right|\Bigr\}+\text{o}_{\hat{\mathbb{P}}}(1).
\end{align*}

From Theorem 1 in \cite{shao1995} it follows that
\begin{align*}
\max_{\substack{1\leq i\leq j\leq n \\ j-i+1\geq nc_n}} \Bigl\{ \sqrt{j-i+1} \left|\overline{U}_i^j\right|\Bigr\}\leq \max_{1\leq i\leq j\leq n} \Bigl\{ \sqrt{j-i+1} \left|\overline{U}_i^j\right|\Bigr\}=\mathcal{O}(\sqrt{\log n}) \quad \text{a.s.}
\end{align*}
In combination with Proposition  \ref{estimate_rate}, this yields \eqref{zuzeigen}. 
\hfill $\Box$

\medskip
\medskip

\subsection{Proof of Theorem \ref{Thm2}} 
Note that by definition of $\hat{K}$, we have 
\begin{align*}
\hat{K}(V_{n},q_n)<K^\ast \quad \Longleftrightarrow \quad \exists \vartheta\in\mathcal{S}_n \text{ with } |J(\vartheta)|<K^\ast \text{ such that } V_{n}(Y,\vartheta)\leq q_n
\end{align*}
[see \cite{frick2014}]. It therefore suffices to show that the probability of the existence of a candidate function $\vartheta\in\mathcal{S}_n$ having less than $K^\ast$ change points and fulfilling $V_{n}(Y,\vartheta)\leq q_n$ converges to $0$.

We will again first prove an analog of the result for the statistic $V_{n,\sigma_\star}$, where the estimator $\hat{\sigma}_\star^2$ is replaced by the long run variance.

\begin{theorem} \label{Thm5}
Under the same assumptions as in Theorem \ref{Thm4}, assume that $\{q_n\}_{n\in\N}$ is a sequence fulfilling \eqref{assumption_concerning_q_n}. Then it follows that
\begin{align*}
\lim_{n\rightarrow\infty} \mathbb{P}\left(\hat{K}(V_{n,\sigma_\star},q_n)<K^\ast\right)=0.
\end{align*}
\end{theorem}

\underline{\textit{Proof of Theorem \ref{Thm5}}~:}
We proceed similarly as in the proof of Theorem 7.10 in \cite{frick2014}. Set
\begin{align*} 
\lambda:=\inf_{0\leq k\leq K^\ast} \left|\tau_{k+1}^\ast-\tau_k^\ast\right|, \quad \beta:=\inf_{1\leq k\leq K^\ast} \left|\theta_k^\ast-\theta_{k-1}^\ast\right|,
\end{align*}
and define the $K^\ast$ disjoint intervals
\begin{align*}
I_i:=\Big[\tau_i^\ast-\frac{\lambda}{2},\tau_i^\ast+\frac{\lambda}{2}\Big), \quad i=1,\ldots,K^\ast. 
\end{align*}
Moreover, define $\theta_i^{+}:=\max\{\theta_{i-1}^\ast,\theta_i^\ast\}$, $\theta_i^{-}:=\min\{\theta_{i-1}^\ast,\theta_i^\ast\}$, $I_i^{+}:=\{t\in I_i\,:\,\vartheta^\ast(t)=\theta_i^{+}\}$, and $I_i^{-}:=\{t\in I_i\,:\,\vartheta^\ast(t)=\theta_i^{-}\}$. Note that $|I_i^{+}|=|I_i^{-}|=\lambda/2$. In particular, since $\{c_n\}_{n\in\N}$ is a null sequence, it holds that $|I_i^{+}|\geq c_n$ and $|I_i^{-}|\geq c_n$ for any $n\in\N$ large enough.

Any candidate function with $K<K^\ast$ change points must be constant on at least one of the disjoint intervals $I_i$. Therefore we get
\begin{align*}
\mathbb{P}\left(\hat{K}(V_{n,\sigma_\star},q_n)<K^\ast\right)
&\leq \sum_{i=1}^{K^\ast} \mathbb{P}\left(\exists \theta\leq \theta_i^{+}-\frac{\beta}{2}\,:\,\frac{1}{\sigma_\star}\sqrt{\frac{n\lambda}{2}}\left|\overline{Y}_{I_i^{+}}-\theta\right|-\sqrt{2\log\frac{2e}{\lambda}}\leq q_n\right)\\
&+\sum_{i=1}^{K^\ast}  \mathbb{P}\left(\exists \theta\geq \theta_i^{-}+\frac{\beta}{2}\,:\,\frac{1}{\sigma_\star}\sqrt{\frac{n\lambda}{2}}\left|\overline{Y}_{I_i^{-}}-\theta\right|-\sqrt{2\log\frac{2e}{\lambda}}\leq q_n\right).
\end{align*}

All of these summands can be dealt with analogously, which is why we will restrict ourselves to the second probability and the case $i=1$. Without loss of generality, assume that $I_1^{-}=[\tau_1^\ast-\lambda/2,\tau_1^\ast)$. It follows easily that the term of interest is bounded from above by 
\begin{align} \label{upper_bound}
\mathbb{P}\left(\frac{1}{\sigma_{\star}}\sqrt{\frac{n\lambda}{2}}\left|\overline{\varepsilon}_{n\tau_1^\ast-\frac{n\lambda}{2}}^{n\tau_1^\ast-1}-\frac{\beta}{2}\right|-\sqrt{2\log\frac{2e}{\lambda}}\leq q_n\right)
+\mathbb{P}\left(\overline{\varepsilon}_{n\tau_1^\ast-\frac{n\lambda}{2}}^{n\tau_1^\ast-1} > \frac{\beta}{2}\right).
\end{align}
Since $\{\varepsilon_i\}_{i\in\Z}$ is mean-ergodic, the second probability in \eqref{upper_bound} converges to $0$. 
Concerning the first probability in \eqref{upper_bound}, note that with exactly the same Gaussian approximation arguments as given in the proof of Theorem \ref{Thm4}, it suffices to show that
\begin{align} \label{suffices_to_show}
\hat{\mathbb{P}}\left(\frac{1}{\sigma_{\star}}\sqrt{\frac{n\lambda}{2}}\left|\overline{U}_{n\tau_1^\ast-\frac{n\lambda}{2}}^{n\tau_1^\ast-1}-\frac{\beta}{2}\right|-\sqrt{2\log\frac{2e}{\lambda}}\leq q_n\right)=\text{o}(1),
\end{align}
where $U_1,\ldots,U_n\sim\mathcal{N}(0,\sigma_\star^2)$ are i.i.d.\ and defined on a richer probability space $(\hat{\Omega},\hat{\mathcal{A}},\hat{\mathbb{P}})$. 
From Theorem 7.10 and Lemma 7.11 in \cite{frick2014}, it follows that the probability in \eqref{suffices_to_show} is upper bounded by 
$ e^{-\frac{1}{64\sigma_\star^2} n\lambda\beta^2+\frac{1}{2}(q_n+\sqrt{2\log\frac{2e}{\lambda}})^2}$.
By assumption \eqref{assumption_concerning_q_n}, this expression vanishes as $n$ tends to $\infty$. 
\hfill $\Box$
\\

\underline{\textit{Proof of Theorem \ref{Thm2}}~:}
With exactly the same arguments as given in the proof of Theorem \ref{Thm5}, it suffices to show that
\begin{align*}
\lim_{n\rightarrow\infty} \mathbb{P}\left(\frac{1}{\hat{\sigma}_\star}\sqrt{\frac{n\lambda}{2}}\left|\overline{\varepsilon}_{n\tau_1^\ast-\frac{n\lambda}{2}}^{n\tau_1^\ast-1}-\frac{\beta}{2}\right|-\sqrt{2\log\frac{2e}{\lambda}}\leq q_n\right)= 0.
\end{align*}
Set 
\begin{align*}
X_n:=\frac{1}{\sigma_\star}\sqrt{\frac{n\lambda}{2}}\left|\overline{\varepsilon}_{n\tau_1^\ast-\frac{n\lambda}{2}}^{n\tau_1^\ast-1}-\frac{\beta}{2}\right| \quad \text{and} \quad Y_n:=\frac{1}{\hat{\sigma}_\star}\sqrt{\frac{n\lambda}{2}}\left|\overline{\varepsilon}_{n\tau_1^\ast-\frac{n\lambda}{2}}^{n\tau_1^\ast-1}-\frac{\beta}{2}\right|.
\end{align*}
Let $\delta>0$ be arbitrary and note that
\begin{align*} 
\mathbb{P}\left(Y_n\leq q_n+\sqrt{2\log\frac{2e}{\lambda}}\right)&=\mathbb{P}\left(Y_n\leq q_n+\sqrt{2\log\frac{2e}{\lambda}}\, , \, \left|\frac{Y_n}{X_n}-1\right|>\delta\right)\\
&+\mathbb{P}\left(Y_n\leq q_n+\sqrt{2\log\frac{2e}{\lambda}}\, , \, \left|\frac{Y_n}{X_n}-1\right|\leq \delta\right). \nonumber
\end{align*}
By Proposition \ref{estimate_rate}, the first probability converges to $0$. Concerning the second probability, it holds that
\begin{align*}
\mathbb{P}\left(Y_n\leq q_n+\sqrt{2\log\frac{2e}{\lambda}}\, , \, \left|\frac{Y_n}{X_n}-1\right|\leq \delta\right) \leq \mathbb{P}\left(\,\left(1-\delta\right)X_n\leq q_n+\sqrt{2\log\frac{2e}{\lambda}}\right).
\end{align*}
With the arguments given in the proof of Theorem \ref{Thm5}, the expression on the right hand side is bounded by $e^{-\frac{1}{64\sigma_\star^2}n\lambda\beta^2+\frac{1}{2}\left((q_n+\sqrt{2\log\frac{2e}{\lambda}})/(1-\delta)\right)^2}+\text{o}(1)$, which converges to $0$ by \eqref{assumption_concerning_q_n} for a fixed $\delta > 0$. 
\hfill $\Box$
\\
\\

\subsection{{Proof of Theorem \ref{Thm3}}}

We will prove a corresponding statement for the statistic $V_{n,\sigma_\star}$ where again the estimator  $\hat \sigma_\star^{2}$ is replaced 
by the long run variance, that is 
\begin{align*}
\lim_{n\rightarrow\infty} \mathbb{P}\left(\sup_{\vartheta\in\mathcal{C}(V_{n,\sigma_\star},q_n)} \max_{\tau^\ast\in J(\vartheta^\ast)} \min_{\tau\in J(\vartheta)} \left|\tau^\ast-\tau\right|>c_n\right)= 0 .
\end{align*}
For a proof of this statement we define the value $\beta$ and the intervals $J_i$, $J_i^{-}$, and $J_i^{+}$ as in the proof of Theorem \ref{Thm5} by replacing $\lambda/2$ with $c_n$ and then using the letter $J$ instead of $I$.
Any candidate function $\vartheta\in\mathcal{S}_n$ with
\begin{align*} 
\max_{\tau^{\ast}\in J(\vartheta^{\ast})}\min_{\tau\in J(\vartheta)} \left|\tau^{\ast}-\tau\right|>c_n
\end{align*}
must be constant on at least one of the disjoint intervals $J_i$. Assume without loss of generality that $J_1^{-}=[\tau_1^\ast-c_n,\tau_1^\ast)$. With the same arguments as given in the proof of Theorem \ref{Thm5}, it suffices to show that
\begin{align} \label{prob}
\lim_{n\rightarrow\infty} \mathbb{P}\left(\frac{1}{\sigma_\star}\sqrt{nc_n}\left|\overline{\varepsilon}_{n\tau_1^\ast-nc_n}^{n\tau_1^\ast-1}-\frac{\beta}{2}\right|-\sqrt{2\log\frac{e}{c_n}}\leq q_n\right)= 0.
\end{align}
Theorem 7.10 and Lemma 7.11 in \cite{frick2014} in combination with the Gaussian approximation arguments from the proof of Theorem \ref{Thm4} yield that the probability in \eqref{prob} is bounded by 
$ e^{-\frac{1}{32\sigma_\star^2}nc_n\beta^2+\frac{1}{2}(q_n+\sqrt{2\log\frac{e}{c_n}})^2}+\text{o}(1)$,
which converges to $0$ by \eqref{assumptions_concerning_qn_2}.

For the proof of Theorem \ref{Thm3} assume again that $J_1^{-}=[\tau_1^\ast-c_n,\tau_1^\ast)$
and   note that the same arguments  show that the assertion follows from the statement
\begin{align*}
\lim_{n\rightarrow\infty} \mathbb{P}\left(\frac{1}{\hat{\sigma}_\star}\sqrt{nc_n}\left|\overline{\varepsilon}_{n\tau_1^\ast-nc_n}^{n\tau_1^\ast-1}-\frac{\beta}{2}\right|-\sqrt{2\log\frac{e}{c_n}}\leq q_n\right)= 0 .
\end{align*}
The proof thus works exactly as the proof of Theorem \ref{Thm2}.
\hfill $\Box$

\subsection{Proof of  Proposition  \ref{estimate_rate}}
Note that it suffices to show that
\begin{align*}
\E\left[\left(\hat{\sigma}_\star^2-\sigma_\star^2\right)^2\right]=\mathcal{O}(n^{-2/3}).
\end{align*}
We proceed as in the proof of Theorem 3 in \cite{wuzhao2007}. Moreover, by assumption (A2), we can apply Lemma 4 and Lemma 5 from \cite{wuzhao2007}. For $2\leq i\leq m_n$ we define
\begin{align*}
W_{ik_n}:=\sum_{j=(i-1)k_n+1}^{ik_n} \varepsilon_j-\sum_{j=(i-2)k_n+1}^{(i-1)k_n} \varepsilon_j
\end{align*}
and 
\begin{align*}
r_{ik_n}:=\sum_{j=(i-1)k_n+1}^{ik_n} \vartheta^{\ast}\left(\frac{j}{n}\right)-\sum_{j=(i-2)k_n+1}^{(i-1)k_n} \vartheta^{\ast}\left(\frac{j}{n}\right).
\end{align*}
By Lemma 4 in \cite{wuzhao2007} it then follows that 
\begin{align*}
\E\left[\left(\hat{\sigma}_\star^2-\sigma_\star^2\right)^2\right]
&=\Big|\Big|\frac{1}{2k_n(m_n-1)}\sum_{i=2}^{m_n} (W_{ik_n}+r_{ik_n})^2-\sigma_\star^2\Big|\Big|_2^2\\
&\leq\Big|\Big|\frac{1}{2k_n(m_n-1)}\sum_{i=2}^{m_n}\Big[ (W_{ik_n}+r_{ik_n})^2-||W_{ik_n}||_2^2\Big]\Big|\Big|_2^2+\mathcal{O}(k_n^{-2}).
\end{align*}
Note that 
\begin{align*}
\Big|\Big|\frac{1}{2k_n(m_n-1)}\sum_{i=2}^{m_n}\Big[ (W_{ik_n}+r_{ik_n})^2-||W_{ik_n}||_2^2\Big]\Big|\Big|_2^2&\leq \frac{1}{4k_n^2(m_n-1)^2}\Big|\Big|\sum_{i=2}^{m_n} \left[(W_{ik_n}+r_{ik_n})^2-W_{ik_n}^2\right]\Big|\Big|_2^2\\
&+\frac{1}{4k_n^2(m_n-1)^2}\Big|\Big|\sum_{i=2}^{m_n} W_{ik_n}^2-(m_n-1)||W_{2k_n}||_2^2\Big|\Big|_2^2.
\end{align*}
By Lemma 5 in \cite{wuzhao2007}, the second term is of the order $\mathcal{O}(m_n^{-1})$. We therefore need to deal with
\begin{align*}
\frac{1}{4k_n^2(m_n-1)^2}\Big|\Big|\sum_{i=2}^{m_n} \left[2W_{ik_n}r_{ik_n}+r_{ik_n}^2\right]\Big|\Big|_2^2.
\end{align*}
Lemma 4 in \cite{wuzhao2007} gives that $||W_{ik_n}||_2=\mathcal{O}(\sqrt{k_n})$ uniformly over $i=2,\ldots,m_n$. Moreover, it holds that $r_{ik_n}=\mathcal{O}(k_n)$ uniformly over $i=2,\ldots,m_n$. 
Note that since $\vartheta^{\ast}$ is piecewise constant with $K^{\ast}<\infty$ change points, the set
\begin{align*}
\left \{i\in\{2,\ldots,m_n\}\,\Big|\,r_{ik_n}\neq 0 \right \}
\end{align*}
contains a finite number of elements, independently of $n\in\N$. Therefore, it follows that
\begin{align*}
\Big|\Big|\sum_{i=2}^{m_n} \left[2 W_{ik_n}r_{ik_n}+r_{ik_n}^2\right]\Big|\Big|_2=\mathcal{O}(k_n^2),
\end{align*}  
which yields that 
\begin{align*}
\frac{1}{4k_n^2(m_n-1)^2} \Big|\Big|\sum_{i=2}^{m_n} \left[2W_{ik_n}r_{ik_n}+r_{ik_n}^2\right]\Big|\Big|_2^2=\mathcal{O}\left(k_n^2 m_n^{-2}\right).
\end{align*}
Since $k_n\asymp n^{1/3}$, the claim follows. \hfill $\Box$

\bigskip

{\bf Acknowledgements} This work has been supported in part by the
Collaborative Research Center ``Statistical modeling of nonlinear
dynamic processes'' (SFB 823, Teilprojekt A1, C1) of the German Research Foundation
(DFG).  The authors would like to thank Wei Biao Wu for helpful discussions regarding physical dependence and 
Thomas Hotz, Axel Munk and Florian Pein for helpful discussions about the algorithmic aspects of SMUCE.

\begin{small}
 \bibliography{smuce}

\begin{thebibliography}{}

\bibitem[Bai and Perron, 1998]{bai1998}
Bai, J. and Perron, P. (1998).
\newblock Estimating and testing linear models with multiple structural
  changes.
\newblock {\em Econometrica}, 66:47--78.

\bibitem[Bai and Perron, 2003]{bai2003}
Bai, J. and Perron, P. (2003).
\newblock Computation and analysis of multiple structural change models.
\newblock {\em J. Appl. Econmetr.}, 18:1--22.

\bibitem[Braun et~al., 2000]{braun2000}
Braun, J.~V., Braun, R.~K., and Muller, H.-G. (2000).
\newblock Multiple changepoint fitting via quasilikelihood, with application to
  {DNA} sequence segmentation.
\newblock {\em Biometrika}, 87:301--314.

\bibitem[Chakar et~al., 2017]{chakar2017}
Chakar, S., Lebarbier, E., L{\'{e}}vy-Leduc, C., and Robin, S. (2017).
\newblock A robust approach for estimating change-points in the mean of an
  {AR}(1) process.
\newblock {\em Bernoulli}, 23:1408--1447.

\bibitem[Cho and Fryzlewicz, 2015]{cho2015}
Cho, H. and Fryzlewicz, P. (2015).
\newblock Multiple-change-point detection for high dimensional time series via
  sparsified binary segmentation.
\newblock {\em J. R. Statist. Soc. B}, 77:475--507.

\bibitem[Ciuperca, 2011]{ciuperca2011}
Ciuperca, G. (2011).
\newblock A general criterion to determine the number of change-points.
\newblock {\em Statist. Probab. Lett.}, 81:1267--1275.

\bibitem[Ciuperca, 2014]{ciuperca2014}
Ciuperca, G. (2014).
\newblock Model selection by {LASSO} methods in a change-point model.
\newblock {\em Stat. Papers}, 55:349--374.

\bibitem[Davis et~al., 2006]{davis2006}
Davis, R.~A., Lee, T. C.~M., and Rodriguez-Yam, G.~A. (2006).
\newblock Structural break estimation for nonstationary time series models.
\newblock {\em J. Am. Statist. Ass.}, 101:223--239.

\bibitem[Dette et~al., 1998]{dette1998}
Dette, H., Munk, A., and Wagner, T. (1998).
\newblock Estimating the variance in nonparametric regression - what is a
  reasonable choice?
\newblock {\em J. Roy. Stat. Soc. B}, 60:751--764.

\bibitem[Frick et~al., 2014]{frick2014}
Frick, K., Munk, A., and Sieling, H. (2014).
\newblock Multiscale change point inference.
\newblock {\em Journal of the Royal Statistical Society, Ser. B},
  76(3):495--580.

\bibitem[Fryzlewicz, 2014]{fryzlewicz2014}
Fryzlewicz, P. (2014).
\newblock Wild binary segmentation for multiple change-point detection.
\newblock {\em Ann. Statist.}, 42:2243--2281.

\bibitem[Hall et~al., 1990]{hall1990}
Hall, P., Kay, J.~W., and Titterington, D.~M. (1990).
\newblock Asymptotically optimal difference-based estimation of variance in
  nonparametric regression.
\newblock {\em Biometrika}, 77:521--528.

\bibitem[Harchaoui and L{\'{e}}vy-Leduc, 2010]{harchaoui2010}
Harchaoui, Z. and L{\'{e}}vy-Leduc, C. (2010).
\newblock Multiple change-point estimation with a total variation penalty.
\newblock {\em J. Am. Statist. Ass.}, 105:1480--1493.

\bibitem[Haynes et~al., 2017]{haynes2017}
Haynes, K., Fearnhead, P., and Eckley, I.~A. (2017).
\newblock A computationally efficient nonparametric approach for changepoint
  detection.
\newblock {\em Stat. Comput.}, 27:1293--1305.

\bibitem[Killick et~al., 2012]{killick2012}
Killick, R., Fearnhead, P., and Eckley, I.~A. (2012).
\newblock Optimal detection of changepoints with a linear computational cost.
\newblock {\em J. Am. Statist. Ass.}, 107:1590--1598.

\bibitem[Kolaczyk and Nowak, 2005]{kolaczyk2005}
Kolaczyk, E.~D. and Nowak, R.~D. (2005).
\newblock Multiscale generalised linear models for nonparametric function
  estimation.
\newblock {\em Biometrika}, 92:119--133.

\bibitem[Korkas and Fryzlewicz, 2017]{korkas2017}
Korkas, K. and Fryzlewicz, P. (2017).
\newblock Multiple change-point detection for non-stationary time series using
  wild binary segmentation.
\newblock {\em Statistica Sinica}, 27:287--311.

\bibitem[Lavielle and Moulines, 2000]{lavielle2000}
Lavielle, M. and Moulines, E. (2000).
\newblock Least-squares estimation of an unknown number of shifts in a time
  series.
\newblock {\em J. Time Series Anal.}, 21:33--59.

\bibitem[Li et~al., 2018]{li2018}
Li, H., Guo, Q., and Munk, A. (2018).
\newblock Multiscale change-point segmentation: Beyond step functions.
\newblock {\em arXiv:1708.03942}.

\bibitem[Li et~al., 2016]{li2016}
Li, H., Munk, A., and Sieling, H. (2016).
\newblock {FDR}-control in multiscale change-point segmentation.
\newblock {\em Electron. J. Statist.}, 10:918--959.

\bibitem[Matteson and James, 2014]{matteson2014}
Matteson, D.~S. and James, N.~A. (2014).
\newblock A nonparametric approach for multiple change point analysis of
  multivariate data.
\newblock {\em J. Am. Statist. Ass.}, 109:334--345.

\bibitem[Pein et~al., 2017a]{rpackage}
Pein, F., Hotz, T., Sieling, H., and Aspelmeier, T. (2017a).
\newblock {\em {stepR}: Multiscale change-point inference}.
\newblock R package version 2.0-1.

\bibitem[Pein et~al., 2017b]{pein2017}
Pein, F., Sieling, H., and Munk, A. (2017b).
\newblock Heterogeneous change point inference.
\newblock {\em Journal of the Royal Statistical Society, Ser. B},
  79(4):1207--1227.

\bibitem[Preuss et~al., 2015]{preuss2015}
Preuss, P., Puchstein, R., and Dette, H. (2015).
\newblock Detection of multiple structural breaks in multivariate time series.
\newblock {\em J. Am. Statist. Ass.}, 110:654--668.

\bibitem[Shao, 1995]{shao1995}
Shao, Q.-M. (1995).
\newblock On a conjecture of r{\'{e}}v{\'{e}}sz.
\newblock {\em Proceedings of the American Mathematical Society},
  123(2):575--582.

\bibitem[Tecuapetla-G{\'{o}}mez, 2015]{rpackage2}
Tecuapetla-G{\'{o}}mez, I. (2015).
\newblock {\em {dbacf}: Autocovariance Estimation via Difference-Based
  Methods}.
\newblock R package version 0.0.0.9000.

\bibitem[Tecuapetla-G{\'{o}}mez and Munk, 2017]{tecuapetla2017}
Tecuapetla-G{\'{o}}mez, I. and Munk, A. (2017).
\newblock Autocovariance estimation in regression with a discontinuous signal
  and m-dependent errors: A difference-based approach.
\newblock {\em Scandinavian Journal of Statistics}, 44(2):346--368.

\bibitem[Wu, 2005]{wu2005}
Wu, B.~W. (2005).
\newblock Nonlinear system theory: Another look at dependence.
\newblock {\em Proceedings of the National Academy of Sciences USA},
  102:14150--14154.

\bibitem[Wu, 2011]{wu2011}
Wu, B.~W. (2011).
\newblock Asymptotic theory for stationary processes.
\newblock {\em Statistics and Its Interface}, 4:207--226.

\bibitem[Wu and Zhao, 2007]{wuzhao2007}
Wu, B.~W. and Zhao, Z. (2007).
\newblock Inference of trends in time series.
\newblock {\em Journal of the Royal Statistical Society, Ser. B},
  69(3):391--410.

\bibitem[Wu and Zhou, 2011]{wuzhou2011}
Wu, B.~W. and Zhou, Z. (2011).
\newblock Gaussian approximations for non-stationary multiple time series.
\newblock {\em Statistica Sinica}, 21:1397--1413.

\bibitem[Wu and Phoumaradi, 2009]{wuphoumaradi2009}
Wu, W.~B. and Phoumaradi, M. (2009).
\newblock Banding sample autocovariance matrices of stationary processes.
\newblock {\em Statistica Sinica}, 19:1755--1768.

\bibitem[Yao, 1988]{yao1988}
Yao, Y. (1988).
\newblock Estimating the number of change-points via {S}chwarz' criterion.
\newblock {\em Statist. Probab. Lett.}, 6:181--189.

\bibitem[Yau and Zhao, 2016]{yau2016}
Yau, C.~Y. and Zhao, Z. (2016).
\newblock Inference for multiple change points in time series via likelihood
  ratio scan statistics.
\newblock {\em J. Roy. Stat. Soc. B}, 78:895--916.

\end{thebibliography}
\end{small}

\end{document}